\documentclass{article}
\usepackage{maamonthly}

\newtheorem{lemma}{Lemma}
\newtheorem{theorem}{Theorem}
\newtheorem{proposition}{Proposition}

\begin{document}

\title{A Short Proof of the Simple Continued Fraction Expansion of $e$}
\author{Henry Cohn}

\centerline{\textbf{\LARGE NOTES}}
\smallskip
\centerline{\large Edited by \textbf{William Adkins}} \bigskip

\maketitle

\section{INTRODUCTION.}

In~\cite{E}, Euler analyzed the Riccati
equation to prove that the number $e$ has the continued fraction
expansion
\[
e = [2,1,2,1,1,4,1,1,6,1,1,8,\dots] = 2 + \frac{1}{\displaystyle 1
+ \frac{\strut 1}{\displaystyle 2 + \frac{\strut 1}{\displaystyle
1 + \frac{\strut 1}{\displaystyle 1 + \frac{\strut
1}{\displaystyle 4 + \cdots \vphantom{\frac{\strut
1}{\displaystyle 1}}}}}}}.
\]
The pattern becomes more elegant if one replaces the initial $2$
with $1,0,1$, which yields the equivalent continued fraction
\begin{equation} \label{expansion} [1,0,1,1,2,1,1,4,1,1,6,1,1,8,1,1,10,1,\dots]
\end{equation} because
\[
1 + \frac{1}{\displaystyle 0 + \frac{\strut 1}{\displaystyle 1 + \cdots}} = 2 + \cdots.
\]

One of the most interesting proofs is due to Hermite; it arose as
a byproduct of his proof of the transcendence of $e$ in~\cite{H}.
(See~\cite{O} for an exposition by Olds.) The purpose of this note
is to present an especially short and direct variant of Hermite's
proof and to explain some of the motivation behind it.

Consider any continued fraction $[a_0,a_1,a_2,\dots]$. Its $i$th
convergent is defined to be the continued fraction
$[a_0,a_1,\dots,a_i]$. One of the most fundamental facts
about\break
continued fractions is that the $i$th convergent equals $p_i/q_i$,
where $p_i$ and $q_i$ can be calculated recursively using
\[
p_n = a_np_{n-1} + p_{n-2}, \qquad q_n = a_nq_{n-1}+q_{n-2},
\]
starting from the initial conditions $p_0 = a_0$, $q_0 = 1$, $p_1
= a_0a_1+1$, and $q_1 = a_1$. For\break
a proof see \cite[p.~130]{HW} or any introduction to continued
fractions.

\newpage

\section{PROOF OF THE EXPANSION.}

Let $[a_0,a_1,a_2,\dots]$ be the continued frac-\break{}tion
(\ref{expansion}). In other words, $a_{3i+1} = 2i$ and $a_{3i} =
a_{3i+2} = 1$, so $p_i$ and $q_i$ are as follows:
\[
\begin{array}{c|ccccccccc}
i & 0 & 1 & 2 & 3 & 4 & 5 & 6 & 7 & 8\\
\hline p_i & 1 & 1 & 2 & 3 & 8 & 11 & 19 & 87 & 106\\
q_i & 1 & 0 & 1 & 1 & 3 & 4 & 7 & 32 & 39
\end{array}
\]
(Note that $q_1 = 0$ so $p_1/q_1$ is undefined, but that will not
be a problem.) Then $p_i$ and $q_i$ satisfy the recurrence
relations
\[
\begin{array}[t]{@{}l@{}c@{}l@{}l@{}c@{}l@{}}
p_{3n} &{}={}& p_{3n-1}+p_{3n-2}, & \qquad q_{3n} &{}={}& q_{3n-1}+q_{3n-2},\\
p_{3n+1} &=& 2np_{3n}+p_{3n-1}, &  \qquad q_{3n+1} &=& 2nq_{3n}+q_{3n-1},\\
p_{3n+2} &=& p_{3n+1}+p_{3n}, &  \qquad q_{3n+2} &=& q_{3n+1}+q_{3n}.
\end{array}
\]
To verify that the continued fraction~(\ref{expansion}) equals
$e$, we must prove that
\[
\lim_{i \to \infty}\frac{p_i}{q_i} = e.
\]
Define the integrals
\begin{eqnarray*} A_n &=& \int_0^1\frac{x^n(x-1)^n}{n!}e^x\,dx,\\
B_n &=& \int_0^1\frac{x^{n+1}(x-1)^n}{n!}e^x\,dx,\\
C_n &=& \int_0^1\frac{x^n(x-1)^{n+1}}{n!}e^x\,dx.
\end{eqnarray*}

\begin{proposition} \label{help}
For $n \ge 0$, $A_n = q_{3n}e-p_{3n}$, $B_n = p_{3n+1}-q_{3n+1}e$,
and $C_n = p_{3n+2}-q_{3n+2}e.$
\end{proposition}

\begin{proof}
In light of the recurrence relations cited earlier, we need only
verify the initial conditions $A_0 = e-1$, $B_0 = 1$, and $C_0 =
2-e$ (which are easy to check) and prove the three recurrence
relations
\begin{eqnarray}
A_n &=& -B_{n-1} - C_{n-1},\label{a}\\
B_n &=& -2nA_n + C_{n-1},\label{b}\\
C_n &=& B_n - A_n.\label{c}
\end{eqnarray}

Of course,~(\ref{c}) is trivial. To prove~(\ref{a}) (i.e.,
$A_n+B_{n-1}+C_{n-1} = 0$) integrate both sides of
\[
\frac{x^n(x-1)^n}{n!}e^x + \frac{x^{n}(x-1)^{n-1}}{(n-1)!}e^x +
\frac{x^{n-1}(x-1)^{n}}{(n-1)!}e^x = \frac{d}{dx}
\left(\frac{x^n(x-1)^n}{n!}e^x\right),
\]
which follows immediately from the product rule for derivatives.
To prove~(\ref{b}) (i.e., $B_n+2nA_n-C_{n-1} = 0$) integrate both
sides of
\[
\vadjust{\newpage} \frac{x^{n+1}(x-1)^n}{n!}e^x + 2n
\frac{x^n(x-1)^n}{n!}e^x - \frac{x^{n-1}(x-1)^n}{(n-1)!}e^x =
\frac{d}{dx} \left(\frac{x^n(x-1)^{n+1}}{n!}e^x\right),
\]
which follows from the product rule and some additional
manipulation. This completes the proof.\qed
\end{proof}

The recurrences~(\ref{a}) and~(\ref{b}) can also be proved by
integration by parts.

\begin{theorem}
$e = [1,0,1,1,2,1,1,4,1,1,6,1,1,8,1,1,10,1,\dots].$
\end{theorem}

\begin{proof}
Clearly $A_n$, $B_n$, and $C_n$ tend to $0$ as $n \to \infty$. It
follows from Proposition~\ref{help} that
\[
\lim_{i \to \infty} q_i e - p_i = 0.
\]
Because $q_i \ge 1$ when $i \ge 2$, we see that
\[
e = \lim_{i \to \infty} \frac{p_{i}}{q_{i}} = [1,0,1,1,2,1,1,4,1,\dots],
\]
as desired.\qed
\end{proof}

\section{MOTIVATION.}

The most surprising aspect of this proof is the integral formulas,
which have no apparent motivation. The difficulty is that the
machinery that led to them has been removed from the final proof.
Hermite wrote down the integrals while studying Pad\'e
approximants, a context in which it is easier to see how one might
think of them.

Pad\'e approximants are certain rational function approximations
to a power series. They are named after Hermite's student Pad\'e,
in whose 1892 thesis~\cite{P} they were studied systematically.
However, for the special case of the exponential function they are
implicit in Hermite's 1873 paper~\cite{H} on the transcendence of
$e$.

We will focus on the power series
\[
e^z = \sum_{k \ge 0} \frac{z^k}{k!},
\]
since it is the relevant one for our purposes. A \textit{Pad\'e
approximant} to $e^z$ of type $(m,n)$ is a rational function
$p(z)/q(z)$ with $p(z)$ and $q(z)$ polynomials such that $\deg
p(z) \le m$, $\deg q(z) \le n$, and
\[
\frac{p(z)}{q(z)} = e^z + O\left(z^{m+n+1}\right)
\]
as $z \to 0$. In other words, the first $m+n+1$ coefficients in
the Taylor series of $p(z)/q(z)$ agree with those for $e^z$. One
cannot expect more agreement than that: $p(z)$ has $m+1$
coefficients and $q(z)$ has $n+1$, so there are $m+n+2$ degrees of
free-dom
in toto, one of which is lost because $p(z)$ and $q(z)$ can be
scaled by the same factor without changing their ratio. Thus, one
expects to be able to match $m+n+1$ coefficients. (Of course this
argument is not rigorous.)

It is easy to see that there can be only one Pad\'e approximant of
type $(m,n)$: if $r(z)/s(z)$ is another, then
\[
\vadjust{\newpage} \frac{p(z)s(z)-q(z)r(z)}{q(z)s(z)} =
\frac{p(z)}{q(z)}-\frac{r(z)}{s(z)} = O\left(z^{m+n+1}\right).
\]
Because $p(z)s(z)-q(z)r(z)$ vanishes to order $m+n+1$ at $z = 0$
but has degree\break
at most $m+n$ (assuming $\deg r(z) \le m$ and $\deg s(z) \le n$),
it must vanish identically. Thus, $p(z)/q(z) = r(z)/s(z)$.

We have no need to deal with the existence of Pad\'e approximants
here, because it will follow from a later argument. One can
compute the Pad\'e approximants of $e^z$ by solving simultaneous
linear equations to determine the coefficients of the numerator
and denominator, after normalizing so the constant term of the
denominator is $1$.

The usefulness of Pad\'e approximants lies in the fact that they
provide a powerful way to approximate a power series. Those of
type $(m,0)$ are simply the partial sums of the series, and the
others are equally natural approximations. If one is interested in
approximating $e$, it makes sense to plug $z = 1$ into the Pad\'e
approximants for $e^z$ and see what happens.

Let $r_{m,n}(z)$ denote the Pad\'e approximant of type $(m,n)$ for
$e^z$. Computing continued fractions reveals that
\begin{eqnarray*}
r_{1,1}(1) &=& [2,1],\\
r_{1,2}(1) &=& [2,1,2],\\
r_{2,1}(1) &=& [2,1,2,1],\\
r_{2,2}(1) &=& [2,1,2,1,1],\\
r_{2,3}(1) &=& [2,1,2,1,1,4],\\
r_{3,2}(1) &=& [2,1,2,1,1,4,1],\\
r_{3,3}(1) &=& [2,1,2,1,1,4,1,1],
\end{eqnarray*}
etc. In other words, when we set $z = 1$, the Pad\'e approximants
of types $(n,n)$, $(n$, $n+1)$, and $(n+1,n)$ appear to give the
convergents to the continued fraction of $e$. There is no reason
to think Hermite approached the problem quite this way, but his
paper~\cite{H} does include some numerical calculations of
approximations to $e$, and it is plausible that his strategy was
informed by patterns he observed in the numbers.

It is not clear how to prove this numerical pattern directly from
the definitions. However, Hermite found an ingenious way to derive
the Pad\'e approximants from integrals, which can be used to prove
it. In fact, it will follow easily from Proposition~\ref{help}.

It is helpful to reformulate the definition as follows. For a
Pad\'e approximant of type $(m,n)$, we are looking for polynomials
$p(z)$ and $q(z)$ of degrees at most $m$ and $n$, respectively,
such that $q(z)e^z - p(z) = O\left(z^{m+n+1}\right)$ as $z \to 0$.
In other words, the function
\[
z \mapsto \frac{q(z)e^z - p(z)}{z^{m+n+1}}
\]
must be holomorphic. (Here that is equivalent to being bounded for
$z$ near $0$. No complex analysis is needed in this article.)

One way to recognize a function as being holomorphic is to write
it as a suitable integral. For example, because
\[
\frac{(z-1)e^z + 1}{z^2} = \int_0^1 xe^{zx} \, dx,
\]
it is clear that $z \mapsto ((z-1)e^z+1)/z^2$ is holomorphic. Of
course, that is unnecessary for such a simple function, but
Hermite\vadjust{\newpage} realized that this technique was quite
powerful. It is not clear how he thought of it, but everyone who
knows calculus has integrated an exponential times a polynomial,
and one can imagine he simply remembered that the answer has
exactly the form we seek.

\begin{lemma} \label{intbyparts}
Let $r(x)$ be a polynomial of degree $k$. Then there are
polynomials $q(z)$ and $p(z)$ of degree at most $k$ such that
\[
\int_0^1 r(x) e^{zx} \, dx = \frac{q(z)e^z-p(z)}{z^{k+1}}.
\]
Specifically,
\begin{eqnarray*}
q(z) &=& r(1)z^k - r'(1)z^{k-1} + r''(1)z^{k-2} - \cdots\\
\nnoalign{and}
p(z) &=& r(0)z^k - r'(0)z^{k-1} + r''(0)z^{k-2} - \cdots.
\end{eqnarray*}
\end{lemma}

\begin{proof}
Integration by parts implies that
\[
\int_0^1 r(x) e^{zx} \, dx = \frac{r(1)e^z-r(0)}{z} -
\frac{1}{z}\int_0^1 r'(x) e^{zx}\, dx,
\]
from which the desired result follows by induction.\qed
\end{proof}

To get a Pad\'e approximant $p(z)/q(z)$ of type $(m,n)$, we want
polynomials $p(z)$\break
and $q(z)$ of degrees $m$ and $n$, respectively, such that
\[
z \mapsto \frac{q(z)e^z-p(z)}{z^{m+n+1}}
\]
is holomorphic. That suggests we should take $k = m+n$ in the
lemma. However, if $r(x)$ is not chosen carefully, then the
degrees of $q(z)$ and $p(z)$ will be too high. To choose $r(x)$,
we examine the explicit formulas
\begin{eqnarray*}
q(z) &=& r(1)z^{m+n} - r'(1)z^{m+n-1} + r''(1)z^{m+n-2} - \cdots\\
\nnoalign{and}
p(z) &=& r(0)z^{m+n} - r'(0)z^{m+n-1} + r''(0)z^{m+n-2} - \cdots.
\end{eqnarray*}
The condition that $\deg q(z) \le n$ simply means $r(x)$ has a
root of order $m$ at $x = 1$, and similarly $\deg p(z) \le m$
means $r(x)$ has a root of order $n$ at $x = 0$. Since $\deg r(x)
= m+n$, our only choice (up to a constant factor) is to take
\[
r(x) = x^n(x-1)^m,
\]
and that polynomial works. Thus,
\[
\int_0^1 x^n(x-1)^m e^{zx} \, dx = \frac{q(z)e^z-p(z)}{z^{m+n+1}},
\]
where $p(z)/q(z)$ is the Pad\'e approximant of type $(m,n)$ to
$e^z$.

\newpage

Setting $z = 1$ recovers the integrals used in the proof of the
continued fraction\break
expansion of $e$, except for the factor of $1/n!$, which simply
makes the answer prettier.\break
Fundamentally, the reason why the factorial appears is that
\[
\frac{d^n}{dx^n}\left(\frac{x^n(x-1)^n}{n!}\right)
\]
has integral coefficients. Note also that up to a change of
variables, this expression is the Rodrigues formula for the
Legendre polynomial of degree $n$ \cite[p.~99]{AAR}.

Natural generalizations of these integrals play a fundamental role
in Hermite's proof of the transcendence of $e$. See \cite[p.~4]{B}
for an especially short version of the proof or chapter~20
of~\cite{S} for a more leisurely account (although the integrals
used there are slightly different from those in this paper).

\begin{affil}
Microsoft Research, Redmond, WA 98052-6399\\
cohn@microsoft.com
\end{affil}

\vfill\eject

\end{document}